\documentclass[12pt]{amsart}
\usepackage{amsmath,amsfonts,amsthm,amsopn,cite,mathrsfs}
\usepackage{epsfig,verbatim}
\usepackage{subfigure}

\newcommand{\tfa}{time-frequency analysis}

\newcommand{\ft}{Fourier transform}
\newcommand{\stft}{short-time Fourier transform}

\newcommand{\tf}{time-frequency}

\newcommand{\fif}{if and only if}
\newcommand{\tfs}{time-frequency shift}

\newcommand{\modsp}{modulation space}
\newcommand{\psdo}{pseudodifferential operator}

\newcommand{\knc}{Kohn--Nirenberg correspondence}

\newtheorem{tm}{Theorem}[section]
\newtheorem{lemma}[tm]{Lemma}
\newtheorem{prop}[tm]{Proposition}
\newtheorem{cor}[tm]{Corollary}

\newcommand{\rems}{\noindent\textsl{REMARKS:}}
\newcommand{\rem}{\noindent\textsl{REMARK:}}

 \theoremstyle{definition}
 \newtheorem{definition}{Definition}

\newcommand{\beqa}{\begin{eqnarray*}}
\newcommand{\eeqa}{\end{eqnarray*}}

\DeclareMathOperator*{\supp}{supp}

\newcommand{\field}[1]{\mathbb{#1}}
\newcommand{\bR}{\field{R}}        
\newcommand{\bN}{\field{N}}        
\newcommand{\bZ}{\field{Z}}        
\newcommand{\bC}{\field{C}}        
\newcommand{\bT}{\field{T}}        %
 \def\cF{\mathcal{F}}              
 \def\cS{\mathcal{S}}

 \def\cB{\mathcal{B}}
 
 \def\cG{\mathcal{G}}
 \def\cM{\mathcal{M}}

 \def\cA{\mathcal{A}}

 \def\cC{\mathcal{C}}
 
 \def\cO{\mathcal{O}}

 \def\cX{\mathcal{X}}
 \def\cZ{\mathcal{Z}}

\def\vgf{V_gf}

\def\rd{\bR^d}

\def\rdd{{\bR^{2d}}}
\def\zdd{{\bZ^{2d}}}

\def\lrd{L^2(\rd)}

\def\zd{\bZ^d}

\def\mvv{M_v^1}

\def\lpm{\ell ^p_m}

\def\Lpm{L ^p_m}

\def\intrd{\int_{\rd}}
\def\intrdd{\int_{\rdd}}

\def\<{\left<}
\def\>{\right>}

\def\inv{^{-1}}

\def\mv1{M_v^1}

\def\Lmpq{L_m^{p,q}}

\def\Mmpq{M_m^{p,q}}
\def\phas{(x,\xi )}

\newcommand{\bull}{$\bullet \quad $}

\newcommand{\mif}{M^{\infty,1}}

\newcommand{\vf}{\varphi}
\newcommand{\bba}{\mathbf{a}}
\newcommand{\bbb}{\mathbf{b}}
\newcommand{\bbc}{\mathbf{c}}

\newcommand{\kn}{Kohn-Nirenberg}
\newcommand{\elv}{\ell ^1_v}
\newcommand{\av}{\cA _v}
\newcommand{\inc}{inverse-closed}
\newcommand{\subm}{submultiplicative}
\begin{document}
\begin{abstract}
We discuss the most common types of weight functions in harmonic
analysis and how they occur in \tfa . As a general rule,
submultiplicative weights characterize algebra properties, moderate
weights  characterize module properties, Gelfand-Raikov-Shilov weights
determine spectral invariance, and Beurling-Domar weights guarantee
the existence of compactly supported test functions. 
\end{abstract}

\title{Weight Functions in Time-Frequency Analysis}
\author{Karlheinz Gr\"ochenig}
\address{Faculty of Mathematics \\
University of Vienna \\
Nordbergstrasse 15 \\
A-1090 Vienna, Austria}
\email{karlheinz.groechenig@univie.ac.at}
\date{}
\thanks{K.~G.~was supported by the Marie-Curie Excellence Grant MEXT-CT 2004-517154}
\maketitle

\section{Introduction}

Weight functions are a very technical topic in \tfa . Many 
different conditions on weights appear in the literature, and their
motivation is sometimes confusing.
 This article offers a survey of the most important classes of  weight
 functions in \tfa .  

Weights are used to quantify growth and decay conditions. For
instance, if $m(t) = (1+|t|)^s$ and $\|f\|_{L^\infty _m} = \sup _{t\in
  \rd } |f(t) | m(t) < \infty$, then $|f(t)| \leq C (1+|t|)^{-s}$. So
if $s>0$, then this condition describes the polynomial decay of $f$ of 
order $s$, whereas if $s<0$, then $f$ grows at most like a polynomial
of degree $s$.  Combining this intuition with $L^p$-spaces, one
obtains the weighted $L^p$-spaces which are defined by the norm
$\|f\|_{\Lpm } = \|f \, m \|_p = \intrd |f(t) | ^p m(t) ^p \, dt $. 

The harmonic analysis of weighted $L^p$-spaces is understood to a
large extent. An important source for convolution relations and
algebra properties is Feichtinger's early paper~\cite{feichtinger79}. 

Weight functions in \tfa\ occur in many problems and contexts: (a) in
the definition of \modsp s where the weight helps to  measure and describe the
\tf\ concentration of a function, (b) in the definition of symbol
classes for \psdo s  where the weight describes the specific form of
smoothness in the Sj\"ostrand class, and (c) in the theory of Gabor
frames and \tf\ expansions  where the weight measures the quality of
\tf\ concentration.

This article is organized as follows: Section~2 contains the
definitions and examples of several classes of weight functions,
Section~3 the main definitions of \tfa . In the subsequent sections we
discuss the class of \subm\ weights, moderate weights, GRS-weights,
subconvolutive weights, and Beurling-Domar weights. In each section,
we recall first the main  definition and then state the characterizing
property in the context of $L^p$-spaces. We then present the main
applications and results about these weights in \tfa . Although we
will not offer complete proofs, we try to sketch those proof ideas
that shed some light  on why a particular weight class arises in
harmonic analysis and \tfa .  




\section{Classes of Weight Functions}

In general a weight function $m$ on $\rd $  is simply a
non-negative function. We will assume without loss of generality that
the weight is continuous.  

In \tfa\ the following types of weight functions occur.

\begin{definition}
  Let $v$ and $m$ be non-negative functions on $\rd $ or on $\zd $. 

(a) A weight $v$ is called \emph{submultivplicative}, if 
\begin{equation}
  \label{eq:2}
  v(x+y) \leq v(x) v(y) \qquad \forall x,y \in \rd \, .
\end{equation}

(b) Given a submultiplicative weight $v$, a non-negative function  $m$ is  called
a \emph{$v$-moderate} weight, if there exists a constant $C> 0$, such that  
\begin{equation}
  \label{eq:3}
  m(x+y ) \leq C v(x) m(y) \qquad \forall x,y \in \rd \, . 
\end{equation}
We denote the set of all $v$-moderate weights by $\cM _v$. If $m$ is
$v$-moderate with respect to some $v$, then we simply call $m$
moderate. 

(c) A non-negative function $v$  is called \emph{subconvolutive}, if
$v\inv \in L^1(\rd )$ and $v\inv \ast v\inv \leq C v\inv $ (as a
pointwise inequality).

(d) A submultiplicative weight $v$ satisfies the \emph{GRS-condition}
(the Gelfand-Raikov-Shilov condition), if
\begin{equation}
  \label{eq:4}
  \lim _{n\to \infty } v(nx) ^{1/n} = 1 \qquad \forall x \in \rd \, .
\end{equation}
Equivalently, $\lim _{n\to \infty } \log v(nx) / n = 0$ for all $x\in
\rd $. 

(e) A submultiplicative weight $v$ satisfies the \emph{Beurling-Domar
  condition} (BD-condition), if 
\begin{equation}
  \label{eq:5}
  \sum _{n=0} ^\infty \frac{\log v(nx)}{n^2} < \infty \qquad \forall x
  \in \rd \, .
\end{equation}

(f) A submultiplicative weight $v$ satisfies the \emph{logarithmic
  integral condition}, if 
\begin{equation}
  \label{eq:6}
  \int _{|x| \geq 1} \frac{\log v(x)}{|x|^{d+1}} \, dx < \infty \, .
\end{equation}
\end{definition}

\rems\ 1. All these conditions (except for (f)) make sense on
arbitrary  locally compact Abelian groups or even on 
 more general groups, but for simplicity we will restrict our
attention to $\rd $,  the corresponding phase-space $\rdd $, and to
$\zd $. 

2. Feichtinger's early paper~\cite{feichtinger79} contains a detailed study of
\subm , moderate, and subconvolutive weights on locally compact groups
with many examples and counter-examples. GRS-weights occur first in
the Russian literature~\cite{gelfandraikov}. The logarithmic
integral condition was found by  Beurling~\cite{beurling38}, it is so
prominent in analysis that Koosis devoted an entire 
monograph to it~\cite{koosis}, the Beurling-Domar condition was
discovered by Domar in~\cite{domar56}. 

\vspace{3 mm}

\textsl{The standard examples.} We consider the following class of
weight functions:
\begin{equation}
  \label{eq:7}
  m(x) = m_{a,b,s,t}(x)  = e^{a |x|^b} \, (1+|x|)^s \, \big( \log
  (e+|x|)\big)^t \,  .
\end{equation}
In particular, this class contains  \\
\bull the \emph{polynomial weights} $m_s
(x) = (1+|x|)^s, s \in \bR$, \\
\bull  the  \emph{exponential weights} $m(x) = e^{a|x|}, a \in \bR$,\\
\bull and the  \emph{subexponential weights} $m(x) = e^{a|x|^b}$ for $0\leq b
< 1$ and $a \in \bR $.  
\vspace{3mm}

 The following classification is taken from ~\cite{feichtinger79}. The
 proof is elementary and omitted. 

\begin{lemma}
  (a) If $a,s,t \geq 0 $ and $0\leq b \leq 1$, then $m=m_{a,b,s,t}$ is
  submultiplicative. 

(b) If $a, s, t \in \bR$ and $0\leq  b \leq 1$, then $m=m_{a,b,s,t}$ is
 moderate. 

(c) If either $0< b <1, a>0, s,t\in \bR$ or $b\in \{0,1\}$  and $s>d$,
then $m_{a,b,s,t}$ is subconvolutive.  

(d) If $a,s,t \geq 0$ and $0\leq b <1$, then  $m_{a,b,s,t}$ satisfies
the GRS-condition, the 
Beurling-Domar condition and the logarithmic integral condition.
   \end{lemma}

\begin{lemma}
  If $v$ satisfies the Beurling-Domar condition, then it satisfies the GRS-condition.
\end{lemma}
\begin{proof}
We argue by contradiction and show that if the GRS-condition is
violated, then the Beurling-Domar condition cannot hold either. 
Since $\log v(x^n)$ is subadditive, a standard lemma~\cite[Lemma~VIII.1.4]{DS88}  implies that
\begin{equation}
  \label{eq:2345}
\lim _{n\to \infty } \frac{\log v(x^n)}{n} = \inf _{n\to \infty }
\frac{\log v(x^n)}{n}   \, .
\end{equation}
 In particular,  this limit exists.
If the GRS-condition is not satisfied, then there exists a subsequence
$n_k$ and $\alpha >0$, such that $\frac{\log v(x^n)}{n} \geq \alpha
>0$ and by \eqref{eq:2345} therefore $\frac{\log v(x^n)}{n} \geq
\frac{\alpha}{2}>0$ for all $n \geq N_0$. Consequently we also have 
$\sum _{n=N_0}^\infty \frac{\log v(x^n)}{n^2} \geq \sum
_{n=N_0}^\infty \frac{\alpha}{2n} = \infty $. So the Beurling-Domar
condition is not satisfied. 
\end{proof}
 \rems\ 1. The weight $v(x) = e^{|x|/\log (e+|x|)}$ is
 submultiplicative; it  satisfies  the
 GRS-condition, but not the  BD-condition. 

2.  It is known that the Beurling-Domar condition and the
logarithmic integral condition are equivalent. 

3. Other examples of weights can be obtained by replacing the
Euclidean norm $| \cdot |$ on $\rd $ by some other norm on $\rd $ and
by restriction to a subspace of $\rd $. As an example of this
procedure we mention $m(x_1, x_2, x_3) = (1 + |x_1|^p + |x_3|^q)^{s}$
on $\bR ^3$ where $1\leq p,q < \infty $ and $s\in \bR $.

\section{The Short-Time Fourier Transform}

If $X = \phas 
\in \rdd $ is  a point in the \tf\ plane, the corresponding \tfs\ is
defined by  
\begin{equation}
  \label{eq:8}
  \pi (X)f(t) = M_\xi T_x f(t) = e^{2\pi i \xi \cdot t} f(t-x) \qquad
  t\in \rd \, .
\end{equation} 

Time-frequency shifts do not commute, instead  they satisfy the canonical
commutation relations
\begin{equation}
\label{eq:48}
M_\xi T_x = e^{2\pi i x\cdot \xi } T_x M_\xi  \qquad \phas \in \rdd
\, .  
\end{equation}

The  transform associated to \tfs s  is the \stft . Let $g$ be a
suitable non-zero  window
function on $\rd $, 
then the \stft\ (STFT) of a function or distribution $f$ is
defined to be 
\begin{eqnarray}
  \label{eq:9}
  \vgf (X) &=& \intrd f(t) \,  \overline{g (t-x)} \, e^{-2\pi i \xi
    \cdot t} \,  dt
  \\
&=& \langle f, \pi (X) g \rangle   = \langle f, M_\xi T_x g
\rangle \, .       \notag
\end{eqnarray}

The STFT is  well defined when  we take $g$ to be in a space of test
functions $\cX $ that is invariant under  \tfs s  and $f$ in the dual
space $\cX '$. For instance, we may take 
$f,g \in \lrd $; or $g\in \cS (\rd )$ (the Schwartz class)  and $f\in
\cS ' (\rd )$ (the tempered distributions). No matter which choice of
test functions we make,  the Gaussian  $\vf (t) =
e^{-\pi t\cdot t}$ will always work as a suitable  window.

The covariance property says that \tfs s are mapped to shifts in the
\tf\ plane, because
\begin{equation}
  \label{eq:56}
  |V_g (\pi (Y)f)(X)| = |V_gf (X-Y)| \, .
\end{equation}

For a detailed discussion of \tfs s and the \stft\ we refer
to~\cite[Ch.~3]{book}. 

\section{Submultiplicative Weights}

Recall that a non-negative function $v$ on $\rd $ is
submultiplicative, if 
\begin{equation*}
  v(x+y) \leq v(x) v(y) \qquad \forall x,y \in \rd \, .
\end{equation*}

\vspace{4 mm}

\begin{center}
\fbox{Submultiplicative weights characterize  (Banach) algebra properties.}  
\end{center}

\vspace{4 mm}

\emph{The standard property in harmonic analysis.} As usual, $L^1_v
(\rd )$ is the Banach space  defined by the norm
$$
\|f\|_{L^1_v} := \intrd |f(t) | v(t) \, dt = \|f \, v \|_1 \, .
$$
Likewise, $\elv (\zd ) $ is defined by the norm
$\|\mathbf{a}\|_{\elv } = \sum _{k\in \zd } |a_k| v(k) $. 
\begin{lemma} \label{weilone}
The space $\ell ^1_v (\zd )$ is a Banach
  algebra under convolution, \fif\    $v$ is \subm\ on $\zd $, and
  $L^1_v (\rd ) $ is  a Banach
  algebra under convolution, \fif\    $v$ is \subm\ on $\rd $.
\end{lemma}

\begin{proof}
We use   $v(l)\leq  v(k) v(l-k)$ and 
estimate in a  straightforward manner:
\begin{eqnarray*}
\|\bba  \ast \bbb \|_{\ell ^1_v} &=& \sum _{l\in \zd } \Big|\sum _{k\in
  \zd }  a_k b_{l-k} \Big| \, v(l) \\ 
&\leq  & \sum _{l \in \zd } \sum _{k \in \zd }| a_k| v(k) \,\,  |
b_{l-k}|   v(l-k ) \\
&=& \Big(\sum _{l\in \zd } |b_l| v(l)  \Big) \, \Big(\sum _{k\in \zd }
|a_k| v(k)  \Big) = \|\bba \|_{\elv }\, \|\bbb \|_{\elv } \, .  
\end{eqnarray*}
To show the converse, let $\delta _m(k) = 1$ for $k=m$ and $\delta
_m(k)=0$ for 
$k\neq m$. Then $\delta _m \ast \delta _n = \delta _{m+n}$ and $\|\delta
_m \|_{\elv } = v(m) $. Thus, if $\elv $ is a Banach algebra, then  
$$
v(m+n) = \|\delta _{m+n} \|_{\elv } = \|\delta _m \ast \delta _n
\|_{\elv } \leq \|\delta _m \|_{\elv } \, \|\delta _n \|_{\elv } =
v(m) v(n) \, ,
$$ 
and so $v$ is \subm . The proof for $L^1_v (\rd ) $ is similar, but
requires approximate identities for the converse. 
\end{proof}

\rems\ 
1. Usually $L^1_v (\rd ) $ is equipped with  the
  involution $f^\ast (x) := \overline{f(-x)}$. It is easy to verify
  that this involution  is   continuous, \fif\ $v(-x) \leq C
  v(x)$. It is therefore convenient, but not absolutely necessary, to
  assume that $v$ is an even function.  If $v $ is even, then  $v(0) \leq v(0)^2
\leq v(x) v(-x) = v(x) ^2$. Thus either  $v(0) = 0$ and  $v\equiv 0$
(because $v(x) \leq v(0) v(x)$), or $v(0) \geq 1$ and then   $v(x) \geq
1$ for all $x$. 

2.  From now on,  we will  assume without loss of generality that
$v$ is an even 
function, because we are mostly interested in involutive Banach
algebras. The treatment of Banach algebras with unbounded involution
is also possible~\cite{palmer1}, but this case  occurs rarely in \tfa .

   \begin{lemma}\label{expo}
If $v$ is submultiplicative (and even), then there exists a constant
$a \geq 0$, such that 
$$
v(x) \leq  e^{a|x|}\, .$$ 
Every submultiplicative weight grows at most exponentially.
   \end{lemma}

   \begin{proof}
  Define $a$ by  $e^a= \sup _{|t| \leq 1 } v(t) $. Since $v$ is
  continuous and $1= v(0) \leq v(x)v(-x)$, we have $a\geq 1$. Given
  $x\in \rd $ arbitrary, choose $n\in \bN $, so that $n-1 < |x|\leq
  n$. Then $|x/n| \leq 1$ and by the submultiplicativity we find that
$$
    v(x) = v(n \, \dfrac{x}{n}) \leq v(\dfrac{x}{n}) ^n \leq  e^{an} \, ,
$$
thus $v$ grow at most exponentially. 
   \end{proof}

In \tfa\ submultiplicative weights occur in the investigation of
twisted convolution,  in the definition of ``good
windows'' and spaces of test functions, and in the construction of algebras  of
pseudodifferential operators. 

\vspace{4 mm}

\subsection{Series of Time-Frequency Shifts and Twisted Convolution} 
Given a lattice $\Lambda = \alpha \zd \times \beta \zd  \subseteq \rdd
$, we  consider  series of 
\tfs s $A= \sum _{k,l \in \zd  } a_{kl} M_{\beta l}
T_{\alpha k} $. To avoid convergence problems,  it is often
convenient to consider absolutely convergent 
series of \tfs s.   This motivates  the following definition. 

\begin{definition}
The linear space  $\cA _v (\alpha ,\beta)$ consists of all series of \tfs s $  A= \sum
_{k,l \in \zd } a_{kl} M_{\beta l} T_{\alpha k}$ with $\mathbf{a} =
(a_{kl} )_{k,l \in \zd  } \in \ell ^1 _v (\zdd )$. 
\end{definition}
Let $\pi (\mathbf{a}) = \sum _{k,l \in \zd  } a_{kl} M_{\beta l}
T_{\alpha k} $ be mapping from coefficients to operators. By definition
$\pi $  maps $\ell ^1_v (\zdd )$ onto $\cA _v (\alpha ,\beta
)$.  It  can be shown that $\pi  $ is
one-to-one~\cite{GL03,rieffel88}. Consequently   $\|A\|_{\cA _v} 
:= \|\mathbf{a}\|_{\ell ^1_v} $ 
is a Banach space norm on $\cA _v (\alpha,\beta )$.

  \begin{lemma}
    If $v$ is submultiplicative, then $\cA _v(\alpha ,\beta )$ is a Banach
    algebra. 
  \end{lemma}
  \begin{proof}
We introduce a new product between two sequences $\mathbf{a}$ and
$\mathbf{b}$: the  \emph{twisted   convolution}  is
defined  by 
\begin{equation}
  \label{eq:54}
  (\bold{a} \,\, \sharp \, \bold{b}  )(k,l) = \sum _{k',l' \in \zd }
  a_{k'l'}   b_{k-k', l-l'} e^{-2\pi i \theta  k' \cdot (l-l')} \, .
\end{equation}
Then a simple computation using the commutation relations~\eqref{eq:48}
shows that the composition of two series of \tfs s corresponds to the
twisted convolution of the coefficient sequences.  Formally, if $A =
\pi (\mathbf{a})$ and $B = \pi (\mathbf{b})$, then  
\begin{equation}
  \label{eq:55}
AB=   \pi (\mathbf{a}) \pi (\mathbf{b}) = \pi
(\mathbf{a} \, \sharp \, \mathbf{b}) \, .
\end{equation}
Twisted convolution is majorized by ordinary convolution via the
obvious  pointwise inequality  
$$
|(\mathbf{a} \, \sharp \,
\mathbf{b})(k,l)| \leq (|\mathbf{a}| \ast |\mathbf{b}|)(k,l),  \qquad  \forall
 k,l \in \zd \, .
$$ 
Therefore  the Banach algebra property of $\ell ^1_v (\zdd )$
implies  that 
$$
\|AB\|_{\cA _v} = \|\pi (\mathbf{a} \, \sharp \, \mathbf{b})\|_{\cA _v
} = \|\bba \, \sharp \, \bbb \|_{\ell^1_v }  \leq  \|\mathbf{a} \|_{\ell
  ^1_v }  \| \mathbf{b}\|_{\ell
  ^1_v } = \|A \|_{\cA _v}  \|B \|_{\cA _v}  \, .
$$
So $\cA _v$ is a Banach algebra. 
  \end{proof}
\rem\ A similar statement holds for series of \tfs s on an arbitrary
lattice. 

\vspace{3 mm}

\subsection{Spaces of Test Functions and ``Good'' Windows}
 In analysis,
test functions are defined by their smoothness and decay conditions,
and the resulting spaces are usually Frechet spaces. In \tfa , the
appropriate spaces of test functions  are defined by properties of the
STFT. 

\begin{definition}
Let $\vf (t) = e^{-\pi t^2}$ be the Gaussian and  $v$ be a
submultiplicative weight on $\rdd $.   The modulation 
space $\mvv (\rd )$  consists of all $f\in \lrd $ for 
  which the norm 
  \begin{equation}
    \label{eq:9a}
    \|f\|_{\mvv } = \intrdd |V_{\vf } f(X)| v(X) \, dX 
  \end{equation}
is finite. 
\end{definition}


The following lemma asserts that  $\mvv $ is always non-trivial.  

\begin{lemma}
  \label{admissible}
If $v$ is submultiplicative, then $\mvv $ is
non-trivial. Specifically, if $\vf (t) = e^{-\pi t^2}$, then $V_{\vf }
\vf \in \mvv $. More generally,  every function of the form $ g = \sum
_{j=1}^n c_j \pi (X_j )\vf  $, where $  c_j \in \bC, X_j \in \rdd $,
belongs to  $\mvv $. 
\end{lemma}

\begin{proof}
A calculation with Gaussian integrals shows that $V_\vf \vf (X) =
2^{-d/2} e^{i\pi x\cdot \xi } \, e^{-\pi X^2
  /2}$~\cite[Lemma~1.5.1]{book}. Since $v$ grows at 
most  exponentially by 
Lemma~\ref{expo}, we find that $V_{\vf } \vf \in L^1_v (\rdd )$ and
thus $\vf \in M^1_v $ for an arbitary \subm\ weight $v$. 

If $g =  \sum
_{j=1}^n c_j \pi (X_j )\vf  $, then $|V_\vf g (X)| \leq  \sum
_{j=1}^n |c_j| |V_{\vf } \vf (X-X_j)|  $ after using~\eqref{eq:56}. We
obtain 
$$
\|g\|_{\mvv  } = \|V_{\vf } g \|_{L^1_v} \leq \sum
_{j=1}^n |c_j|  \, \|V_{\vf } \vf (\cdot -X_j)\|_{L^1_v} \leq \sum
_{j=1}^n |c_j|\,  v(X_j) \,  \|V_{\vf } \vf \|_{L^1_v} < \infty \, .
$$
So $g\in \mvv $. 
\end{proof}

\subsection{ Symbol Classes  for Pseudodifferential Operators}
 We consider \psdo s in the
Kohn-Nirenberg correspondence. Given a \emph{symbol} $\sigma $ on
$\rdd $, the \psdo\ $K_\sigma $ is defined by 
\begin{equation}
  \label{eq:10}
  K_\sigma f(t) = \intrd \sigma \phas \hat{f} (\xi ) e^{2\pi i
    x\cdot \xi } \, d\xi \, ,
\end{equation}
whenever the integral makes sense. In the tradition of PDE, this operator is
written as $\sigma (x,D)$. Using a suitable duality, 
 the \kn\ correspondence can be defined for an  arbitrary
tempered distribution $\sigma \in \cS ' (\rdd )$, and for even more
general distribution classes. 

After the fundamental papers of J.~Sj\"ostrand~\cite{Sjo94,Sjo95}, the
following symbol classes  have gained some prominence in the
investigation of \psdo s.  

\begin{definition}
  Let $v$ be a submultiplicative (and even)  weight function on $\rdd
  $ and $\Phi (X) = e^{-\pi X^2}$. The \emph{weighted Sj\"ostrand class} $\mif _v (\rdd )$ is
  defined by the norm
  \begin{equation}
    \label{eq:11}
    \|\sigma \|_{\mif _v} = \intrdd \, \sup _{X \in \rdd } |V_\Phi
    \sigma (X,\Xi )| \, v(\Xi ) \, d\Xi \, .
  \end{equation}
Then $\mif _v$ is a subspace of $\cS ' (\rdd )$, consisting of bounded
function that are locally in the Fourier algebra. 
\end{definition}

Composition of \psdo s defines a product on the level of symbols as
follows
\begin{equation}
  \label{eq:12}
K_{\sigma \circ \tau }=   K_\sigma K_\tau   \, .
\end{equation}
If $\sigma , \tau \in \cS (\rdd )$, then  the product is given by the
explicit formula~\cite{hormander3} 
\begin{equation}
  \label{eq:13}
  (\sigma \circ \tau )\phas  = \intrdd \sigma (y,\xi  +\eta ) \tau (x+y,
  \xi ) e^{-2\pi i y\cdot \eta} \, dyd\eta\,  .
\end{equation}

\rem\ There are many other calculi of \psdo s, the most important one
is  the Weyl calculus. The results discussed in this survey are
independent of the used calculus. We use the \knc , because  it is
universal and  can be formulated on arbitrary
locally compact abelian groups. 

\vspace{3 mm}

One of the key properties of the Sj\"ostrand classes is the algebra
property. The following theorem was proved in \cite{Sjo94} for the
unweighted case $v\equiv 1$, a different proof was given
in~\cite{toft01}, weighted versions and two  genuine \tf\ proofs were
obtained in \cite{grocomp,gro06} (for the Weyl calculus), the
extension to LCA groups is contained in \cite{GS06}.

  \begin{tm}\label{sjowei}
  If $v$ is submultiplicative  on $\rdd  $, then $\mif _v
  (\rdd )$
  is a Banach algebra with respect to the  
 product $\, \circ \, $.   
 Furthermore,  
$$  
    \|\sigma \, \circ \,   \tau \|_{\mif _v } \leq C_\Phi \, \|\sigma
    \|_{\mif _v }  
\, \|\tau \|_{\mif _v} \quad \quad \forall  \sigma ,\tau \in \mif _v
(\rdd ) \, . 
  $$
\end{tm}

\begin{proof} We sketch the proof as it is given in
  ~\cite{grocomp}. We   define the \emph{``grand symbol''} of $K_\sigma $ by 
  \begin{equation}
    \label{eq:57}
    G(\sigma ) (\Xi ) = \sup _{X\in \rdd } |V_\Phi \sigma (X, \Xi )|
    \, .
  \end{equation}
Then $\|\sigma \|_{\mif _v} =  \|G(\sigma ) \|_{L^1_v} $ by definition
of $\mif _v $. The technical and difficult part of the proof is to
show the pointwise inequality 
\begin{equation}
  \label{eq:58}
  G(\sigma \circ \tau ) (\Xi ) \leq \Big(G(\sigma ) \ast G(\tau ) \ast H\Big)
  (\Xi ) \qquad \forall \Xi \, \in \rdd \, ,
\end{equation}
where $\ast $ is the ordinary convolution on $\rdd $ and $H$ is a
positive function depending only on $\Phi $. Now the algebra property
of $\mif _v $ follows from the algebra property of $L^1_v (\rdd )$
stated in 
Lemma~\ref{weilone}, because
\begin{eqnarray*}
  \|\sigma \circ \tau \|_{\mif _v}&=& \|G(\sigma \circ \tau )
  \|_{L^1_v } \\
&\leq & \| G(\sigma ) \ast G(\tau ) \ast H \|_{L^1_v} \\
&\leq & \| G(\sigma )  \|_{L^1_v}\, \| G(\tau )  \|_{L^1_v}\, \|  H
\|_{L^1_v} \\
&=& C \|\sigma \|_{\mif _v} \, \|\tau \|_{\mif _v} \, .
\end{eqnarray*}
\end{proof}

\rem\ Let us emphasize that the above result works for arbitrary
submultiplicative weights,  including  exponential weights.

\section{Moderate Weights}

Moderate weights comprise a more general class of weight
functions and are always associated to a \subm\ weight $v$. Precisely,
a non-negative function is called $v$-moderate, if 
$$
  m(x+y ) \leq C v(x) m(y) \qquad \forall x,y \in \rd \, . 
$$
We write $\cM _v$ for the class of all $v$-moderate weight functions. 
From this definition follows that 
$$\frac{1}{C v(x)} \leq m(x) \leq C
v(x)\, ,$$
 and so a moderate weight can grow only as fast as the
associated \subm\ weight $v$. Furthermore, if $K\subseteq \rd $ is
compact, then  $\sup _{t\in K} m(x+t) \leq C \sup _{t\in K} v(t)\,
m(x) = C' m(x)$ and likewise $ \inf _{t\in K} m(x+t) \geq C'' m(x)$. 

Alternatively, a non-negative function is moderate,  
 \fif\
$$
\sup _{y\in \rd } \frac{m(x+y)}{m(y)} = C(x) < \infty \qquad \qquad
\forall x\in \rd \, .
$$
This definition does not make reference to a \subm\ weight $v$. 
A related condition occurs in~\cite[18.4.2]{hormander}.

\vspace{3 mm}

\begin{center}
  \fbox{Moderate weights arise in ``module properties''.}
\end{center}

\vspace{3mm}

Given a weight $m$, we define the weighted $L^p$-space $\Lpm $ 
by the norm
\begin{equation}
  \label{eq:13}
  \|f\|_{\Lpm } := \|f \, m \|_p =  \Big( \intrd |f(t) |^p \, m(t)^p
  \, dt \Big)^{1/p} \,. 
\end{equation}
If $p=\infty $, then $f\in L^\infty _m$ means that $|f(t)|
\leq \|f\|_{L^\infty _m} \, m(t)\inv $.  Likewise $\lpm (\zd ) $ is
defined by the norm $\|\mathbf{a}\|_{\lpm } = \|\mathbf{a}\, m \|_{\ell
  ^p}$.  

The moderateness of a weight  is exactly the condition required for
convolution estimates in the style of Young's inequality~\cite{feichtinger79}. 

\begin{lemma} \label{young}
  Assume that $v$ is \subm\ on $\zd $. Then the following are
  equivalent: 

(i)  $m$ is $v$-moderate. 
 
(ii) $\lpm (\zd )$ is invariant under the shift $T_k$ and the operator norm
satisfies
$$\| T_k\|_{\lpm \to \lpm } \leq C v(k) \qquad \text{ for } \,  k\in \zd \, .$$

(iii) The convolution relation  $\ell^1_v \ast \lpm \subseteq \lpm $
holds  in the sense of the  norm estimate
  \begin{equation}
    \label{eq:14}
    \|\bba \ast \bbb \|_{\lpm} \leq C \|\bba \|_{\elv } \, \|\bbb \|_{\lpm } \, .
  \end{equation}
Similarly, on $\rd $  the convolution relation $L^1_v \ast \Lpm \subseteq
\Lpm $ holds  \fif\ $m\in \cM _v$. 
\end{lemma}

\begin{proof}
The argument is similar to Lemma~\ref{weilone}.  
 See~\cite{feichtinger79} for a detailed proof.  
\end{proof}

Whereas $L^1_v $ is a Banach \emph{algebra},  $\Lpm $ is only
a Banach \emph{space} in general. Lemma~\ref{young} states that $L^1_v$ acts
continuously on $\Lpm $ by convolution. In algebraic terminology, $\Lpm
$ is an $L^1_v$-convolution-\emph{module}. The lemma is the prototype
and explains why and where moderate weights occur. 

\vspace{3 mm}

\subsection{Modulation Spaces}
 In \tfa\ moderate weights occur in the
definition of general \modsp 
s. These spaces are defined in terms of a function space norm applied
to the \stft\ $\vgf $. The idea is to   measure the \tf\ concentration of a
function or distribution. There are several equivalent definitions. We use the most
general definition as discussed in the general theory of coorbit
spaces~\cite{fg89jfa}.

\begin{definition} \label{moddef}
Assume that $v$ is \subm\  and choose $g\in \mvv , g\neq 0$.  Let
$(\mvv )^{\widetilde{}}$ be the space of all 
conjugate linear functionals on $\mvv $. Let $m\in \cM _v$ and  $1\leq
p,q \leq 
\infty $. Then the \modsp\ $\Mmpq (\rd ) $ consists of all $f\in
(\mvv ) ^{\widetilde{}} \, (\rd )$, such that $V_g f \in L^{p,q}_m
(\rdd )$, and
the 
norm is 
\begin{equation}
  \label{eq:17}
  \|f\|_{\Mmpq } = \|\vgf \|_{\Lmpq } =  \Big( \intrd \Big( \intrd |\vgf (x,\xi )|^p \,
  m(x,\xi )^p \, dx \Big)^{q/p} d\xi \Big)^{1/q} \, .
\end{equation}
If $p= \infty $ or $q=\infty $, then we use the supremum norm.   
\end{definition}

We now  state the main properties; these  are essential for a meaningful
theory of these function and distribution spaces~\cite[Ch.~11]{book}. 

\begin{tm} \label{propmod}
  Let $m\in \cM _v$, $1\leq p,q \leq
\infty $. 

(a) Then the \modsp\ $\Mmpq (\rd ) $ is a Banach space. 

(b) If $h\in \mvv , h\neq 0$, then 
\begin{equation}
  \label{eq:18}
  \|V_h f \|_{\Lmpq } \asymp \|V_gf \|_{\Lmpq} \, .
\end{equation}
Thus the definition of $\Mmpq $ is independent of the window $g$, and
different window in $\mvv $ yield equivalent norms on $\Mmpq $. 

(c) Invariance under \tfs s: If $X = (x,\xi ) \in \rdd $, then 
\begin{equation}
  \label{eq:19}
  \|\pi (X) f \|_{\Mmpq } \leq C v(X) \|f\|_{\Mmpq } \, .
\end{equation}
\end{tm}

Let us briefly sketch why the weight $m$ in Definition~\ref{moddef} and
Theorem~\ref{propmod} has to be moderate. 
We use the covariance property of the STFT~\eqref{eq:56} in the form
$\|V_g (\pi (X)f)(Y)| = |V_gf(Y-X)| = T_X |V_gf |(Y)$.  Consequently,
if $m$ is $v$-moderate, then by Lemma~\ref{young}  
$$
\|\pi (X) f \|_{M^{p,p}_m } = \|V_g(\pi (X) f )\|_{\Lpm } = \| T_X V_g f
\|_{\Lpm } \leq C v(X) \|f\|_{M^{p,p}_m } \, . 
$$
Thus  the translation invariance of $\Lpm $ implies the
invariance of $M^{p,p}_m $ under \tfs s. 

The norm equivalence~\eqref{eq:18} is based on a fundamental pointwise
inequaliy for STFTs  \cite[Lemma~11.3.2]{book} 
\begin{equation}
  \label{eq:58}
|V_h f (X)  | \leq \Big( |V_g f | \ast |V_hg| \Big) (X) \qquad \qquad
X\in \rdd \, .    
\end{equation}
If $V_g f \in \Lpm $ and $V_h g \in L^1_v $, then $V_h f \in \Lpm $ by
Lemma~\ref{young}, and once again $m$ needs to be $v$-moderate.

For a meaningful definition of the \modsp s we need (a) a reasonable space
of \emph{test functions} that is invariant under \tfs s, and (b) a
corresponding space of \emph{distributions}. Then the STFT is well
defined and $\Mmpq $ consists of all distributions such that $V_g f
\in \Lmpq $. 
Definition~\ref{moddef} is formulated with the spaces of test
functions and distributions that are intrinsic to \tfa . Specifically,
the appropriate spaces of test functions in \tfa\ are the spaces of
"good windows" $\mvv $,  and the appropriate spaces of distributions are
the dual spaces $(\mvv )^{\widetilde{}} $. This approach works even for
the exponential weights $v(X) = e^{a|X|}$. 

For reasons of convenience, many authors have treated \modsp s for a 
more restrictive class of weight functions, so that standard concepts
can be used.

\begin{prop}
  If $v$ grows at most polynomially, i.e., $v(X) \leq C (1+|X|)^N, X\in
  \rdd $ for some constants $C,N \geq 0$, and $m\in \cM _v$,  then
  $\Mmpq $ is a subspace of the tempered distributions, and so   
  \begin{equation}
    \label{eq:21}
    \Mmpq = \{ f \in \cS ' (\rd ) : \vgf \in \Lmpq \} \, ,
  \end{equation}
whenever $g \in \cS , g\neq 0$.
\end{prop}

Several alternative spaces of test functions have been proposed after
the coorbit space approach in~\cite{fg89jfa}:

1. Let $\vf (t)  = e^{-\pi t^2} $ be the Gaussian on $\rd$, then 
\begin{equation}
  \label{eq:62}
  \cS _{\cC} = \{ g \in \lrd : g = \intrdd F(X) \pi (X) \vf \,
dX , \,\, \mathrm{supp}\, F \, \, \text{is compact}\}
\end{equation}
is a suitable space of test functions that
works for all exponential weights~\cite[Chpt.~11.4]{book}. 

2. In a similar spirit, a discrete version of \eqref{eq:62}, namely
$$ \cS _{\cC , F}= \{ g: g  = \sum
_{j=1}^n c_j \pi (X_j )\vf  \,\, \text{for }   c_j \in \bC, X_j \in
\rdd \} \, ,
$$
 was proposed  in~\cite{CG05} as a suitable window class. Then $\Mmpq $
 can also be  defined as the norm 
completion of $\cS _{\cC , F} $ with respect to the $\Mmpq $-norm. 

3. Another choice is the Gelfand-Shilov space $S^{1/2,1/2}$, which is  used
in~\cite{CPRT}.

For any \subm\ weight $v$, we have the embeddings  $\cS _{\cC , F}
\subseteq S^{1/2, 1/2} \subseteq \mvv  $ and  $\cS _{\cC }
\subseteq S^{1/2, 1/2} \subseteq \mvv  $, and
by~\cite[Prop.~11.4.2]{book} each of $\cS _{\cC}, \cS _{\cC , F}$ and
$S^{1/2 , 1/2}$ is dense in $\mvv $. Consequently, these spaces of
  "special windows'' are universal and work for an  arbitrary \subm\
  weight function $v$.


Beside the natural distribution spaces $M^\infty _{1/v} = (\mvv
)^{\widetilde{}} $ several other spaces have been studied recently. 
In some applications  $m$ is allowed to grow  faster than
polynomially, then  it is
necessary to leave the realm of the Schwartz class  and tempered
distributions. If  $m(X) = \cO (e^{a|X|^b}), 
0\leq b< 1$, then   the  correct
distribution space is the space of  
ultradistributions of Bj\"orck~\cite{bjork66} and
Komatsu~\cite{komatsu73}.  For this reason,  the 
$\Mmpq $ were renamed ultra-modulation spaces by Pilipovic and
Teofanov~\cite{PT04}. If $m$ grows exponentially, e.g., $m(X) = e^{a|X|}$, then
it can be shown that $\Mmpq $ is contained  in the Gelfand-Shilov
space $(S^{1/2,1/2})'$~\cite{CPRT}.
Or to say it differently, all \modsp s with a moderate weight function
are contained in $(S^{1/2,1/2})'$.







With an appropriate distribution space at hand,  \modsp s can be defined
for arbitrary moderate weight functions $m$. In particular \emph{no restriction
  needs to be imposed on the growth of the  moderate weight $m$.}

\vspace{3 mm}

\subsection{Twisted Convolution of Weighted $\ell ^p$-Spaces}

\begin{prop}
  If $m$ is $v$-moderate, then $\ell ^1_v \, \sharp \,  \lpm \subseteq \lpm $
  with  the norm estimate
  \begin{equation}
    \label{eq:14a}
    \| \bba \, \sharp \,  \bbb \|_{\lpm} \leq C \|\bba \|_{\ell ^1_v }
    \, \|\bbb \|_{\lpm } \, . 
  \end{equation}
\end{prop}

\begin{proof}
  Since $|(\bba \, \sharp \, \bbb ) (k,l) | \leq (|\bba | \ast |\bbb |)(k,l)$,
  \eqref{eq:14a} follows from  Theorem~\ref{young} (Young's inequality). 
\end{proof}

\vspace{3 mm}

\subsection{Twisted Product between Modulation Spaces}

A "module property" of \modsp s with respect to the twisted product
can be formulated as follows. 

  \begin{tm}\label{compwei}
  If $v$ is submultiplicative  on $\rdd  $ and $m\in \cM
  _v$,  then $  M^{\infty , p} _m$ is an $\mif _v$-module with respect
  to $\, \circ \, $. This means that the   Young-type  inequality
$$  
    \|\sigma \, \circ \,   \tau \|_{M^{\infty , p }_m } \leq C \, \|\sigma
    \|_{\mif _v }  
\, \|\tau \|_{M^{\infty ,p} _m}\, 
  $$
holds  for all $\sigma \in \mvv $ and all 
  $\tau \in M^{\infty , p } _m$.
\end{tm}

\begin{proof}
The proof  follows from estimate~\eqref{eq:58} for the
  grand symbols and Young's Theorem~\ref{young}. 
\end{proof}


\section{GRS-Weights}

A \subm\ weight $v$ satisfies the Gelfand-Raikov-Shilov condition  
(GRS), if 
$$
\lim _{n\to \infty } v(nx) ^{1/n} = 1 \qquad \forall x\in \rd \, .
$$
The subexponential weight $e^{a |x|^b}$ for $a>0$ and $0\leq b < 1$
  satisfies the GRS-condition, but the exponential weight $e^{a|x|}$
  violates the GRS-condition. Intuitively, the GRS-condition
  describes  the   subexponential growth
  of a weight in   precise technical terms    and excludes all forms
  of exponential growth.  

\vspace{3 mm}

\begin{center}
  \fbox{GRS-weights characterize spectral invariance. }
\end{center}

\vspace{3 mm}

\subsection{Wiener's Lemma}
To prepare the background for the GRS-condition,   we first recall the original 
version of Wiener's lemma for absolutely convergent Fourier series. 

\begin{tm}
  Assume that $f$ is an absolutely converging Fourier series such that
  $f(t) \neq 0$ for all $t \in \bT$, then $1/f$ is also an absolutely
  convergent Fourier series. 
\end{tm}

For  weighted versions of Wiener's Lemma, we define  
\begin{equation*}
  \cA _v (\bT ^d) = \{ f : f(t) = \sum _{k\in \zd} a_k e^{2\pi i k\cdot
    t}, \,  \mathbf{a} = (a_k) \in \ell ^1_v (\zd )\} \, 
\end{equation*}
to be the  space of weighted absolutely convergent Fourier series. 
We equip $\cA _v (\bT ^d)$  with the  norm $\|f\|_{\cA _v} =
\|\textbf{a}\|_{\elv } $, then 
$\cA _v$ is a Banach algebra with respect to pointwise multiplication. 

The weighted version of Wiener's Lemma requires the
GRS-condition and is taken from~\cite{gelfandraikov}.

\begin{tm}
  Assume that $v$ is a submultiplicative weight on $\zd $ satisfying
  the GRS-condition. If $f \in \cA _v(\bT ^d)$ and $f(t) \neq 0$ for
  all $t \in \bT ^d $, then $1/f$ is also in $\cA _v (\bT ^d)$.  
\end{tm}

\subsection{Inverse-Closedness}
 Though not immediate, Wiener's Lemma
should be understood as
 a statement about the relation between two 
Banach algebras. We first define the abstract concept.

\begin{definition}\label{invcl}
  Let $\cA \subseteq \cB $ be two Banach algebras with a common
  identity $e$. We say that $\cA $ is \inc\ in $\cB $, if 
  \begin{equation}
    \label{eq:39}
    a\in \cA \, \, \mathrm{and} \, \, a\inv \in \cB \,\, \Rightarrow a\inv
    \in \cA \, .
  \end{equation}
\end{definition}
In other words, ``the invertibility in the big algebra implies the
invertibility  in the small algebra.'' 

In the case of Wiener's Lemma we take $\cA = \cA _v (\bT ^d)$ and $\cB
= C(\bT ^d)$ or $\cB= L^\infty (\bT ^d)$ (all with pointwise
multiplication). If $f\in \cA _v$ does not vanish anywhere, then by
continuity $\inf _{t\in \bT ^d} |f(t)| >0$ and $f$ is invertible in
$C(\bT ^d)$ (or in $L^\infty(\bT ^d)$). Wiener's Lemma says that $1/f $
must already be in the small algebra $ \cA _v (\bT ^d)$ of weighted  absolutely
convergent Fourier series. So Wiener's Lemma can be recast by saying
that $\cA _v (\bT ^d) $ is inverse-closed in $C(\bT ^d)$. 

 Naimark~\cite{naimark72} turned Wiener's Lemma into a definition and calls
 a nested pair of Banach algebras with common identity a  \emph{Wiener
  pair}, if $\cA $ is \inc\ in $\cB $.

Inverse-closedness can be understood as a spectral property. Let
$\sigma _{\cA } (a) = \{\lambda \in \bC : a-\lambda e \, \text{ is not
invertible in } \, \cA \}$ denote the spectrum of an element $a\in \cA
$. If $\cA = C(\bT ^d)$, then $\sigma _{C(\bT ^d)}(f) = \mathrm{ran}\, f
= f(\bT ^d)$. 

 The following statement is a simple reformulation of Definition~\ref{invcl}.

\begin{lemma}[Spectral invariance, spectral permanence]
  Let $\cA \subseteq \cB $ be two Banach algebras with a common
  identity $e$. Then  $\cA $ is \inc\ in $\cB$, \fif\ 
  \begin{equation}
    \label{eq:40}
    \sigma _{\cA } (a) =     \sigma _{\cB } (a) \qquad \forall a\in
    \cA \, . 
  \end{equation}
\end{lemma}

We now have enough abstract background to reformulate  Wiener's Lemma
and understand the full meaning of  the GRS-condition.

\begin{tm} \label{grr}
  The spectral identity $\sigma _{\av (\bT ^d)} (f) = \sigma _{C(\bT
    ^d)} (f) = \, \mathrm{ran}\, f$ holds for all $f\in \av (\bT ^d)$
  \fif\ $v$ satisfies   the GRS-condition.
\end{tm}

\emph{Counter-Example.} Whereas the sufficiency of the GRS-condition
is quite subtle 
and is based on Gelfand theory for commutative Banach algebras and
complex analysis, the necessity is easy to understand by a
counter-example. The following  reveals the essence
of  the GRS-condition.

Assume that $v$ violates the GRS-condition. Then there exists a $k\in
\zd $ such that $\lim _{n\to \infty } v(nk) ^{1/n} = a = e^\alpha >
1$, and so $v(nk) \geq e^{\alpha n/2}$ for $n\geq n_0$. Thus $v$ grows
exponentially in some direction. 

Let $0< \delta < \alpha /2$ and  set $f(t) = 1 - e^{-\delta } e^{2\pi i
  k\cdot t}$. Since $f$ is a trigonometric polynomial, we have $f\in
\cA _v(\bT ^d)$ and clearly $|f(t)| \geq 1 - e^{-\delta } >  0 $ for
all $t$. The inverse of 
$f$ is the  geometric series 
\begin{equation}
  \label{eq:50}
  \frac{1}{f(t)} = \frac{1}{1 - e^{-\delta } e^{2\pi i
  k\cdot t}} = \sum _{n=0}^\infty e^{-n\delta } e^{2\pi i nk\cdot t}
\, .
\end{equation}
So $1/f$ is an absolutely convergent Fourier series, but 
$$
\|1/f\|_{\cA
  _v} = \sum _{n=0} ^\infty e^{-n\delta } v(nk) \geq \sum _{n=n_0}
^\infty e^{-n\delta } e^{n\alpha /2 } = \infty \, ,
$$
 and so $1/f \not \in
\cA _v(\bT ^d)$.  

\vspace{3 mm}

\subsection{Convolution Operators on $L^1_v$}
 
Wiener's Lemma can be interpreted as a statement about convolution 
operators. Let $C_\mathbf{a}$ be the convolution operator defined by
$C_\mathbf{a}\mathbf{c} = \bba \ast \mathbf{c} $ for $\mathbf{a} \in \ell ^1(\zd
)$ and $ \mathbf{c} \in
\ell ^2(\zd )$. We  identify $\ell ^1 _v  (\zd )$ with the
subalgebra $\mathrm{Op}\, (\ell ^1_v) := \{ C_\mathbf{a} : \mathbf{a} \in \ell ^1_v
(\zd )\}$ of  the C$^*$-algebra $\cB (\ell ^2(\zd ))$ of all bounded operators
on $\ell ^2(\zd )$.  Then Theorem~\ref{grr} is
equivalent to the following.

\begin{tm}\label{wlconv}
  The spectral invariance $\sigma  _{\ell ^1_v } (\mathbf{a}) =
  \sigma_{\cB (\ell ^2)} (C_{\mathbf{a}}) = \mathrm{ran}\,
    \widehat{\mathbf{a}} $ holds for all $\mathbf{a} \in \ell ^1_v (\zd
    )$, \fif\ $v$ satisfies the GRS-condition.
\end{tm}
To see how  this statement follows from Theorem~\ref{grr}, we take
Fourier series $\widehat{\bba } (\xi ) = \sum _{k\in \zd } a_k e^{2\pi
  i k\cdot \xi } $. Then  $(C_{\mathbf{a}}\mathbf{c})\, \widehat{}
\, = \widehat{\mathbf{a}}\, \widehat{\mathbf{c}}$ and  $C_\mathbf{a} $
is unitarily equivalent to the multiplication operator by
$\widehat{\mathbf{a}}$, which has the spectrum $\mathrm{ran} \,
\widehat{\mathbf{a}}$. 

 Replacing the group $\zd $ by some possibly non-commutative locally
compact group $G$, we may ask whether and for which groups  a version of
Theorem~\ref{grr}   still holds. 
The most general result we know of characterizes again GRS-weights.

\begin{tm}[\cite{FGL06}]
  \label{pg}
Let $G$ be a compactly generated group of polynomial growth and $v$ a
submultiplicative weight on $G$. Then the  spectral identity $\sigma
_{L^1_v(G)} (f) = \sigma _{\cB (L^2)} (C_f)$ holds for all $f\in
L^1_v (G)$,  \fif\ $v$ satisfies
  the GRS-condition on $G$, which is $\lim _{n\to \infty }
  v(x^n)^{1/n} = 1 $   for all $x\in G$. 
\end{tm}

\subsection{Wiener's Lemma for Twisted Convolution} 

The version of Wiener's Lemma in the form of Theorem~\ref{wlconv}
suggests that we study the spectrum of the twisted convolution
operator $C_\bba \bbc = \bba \, \sharp \bbc $. Although Fourier series
are no longer available for the non-commutative convolution $\, \sharp
\, $, Wiener's Lemma holds also for twisted convolution. Again, we
identify the algebra $\ell ^1_v (\zdd )$ with a subalgebra of $\cB
(\ell ^2)$ via the isomorphism $\bba \to C_\bba $.

\begin{tm}
  \label{invtwist}
The algebra  $\ell ^1_v (\zdd )$ is
inverse-closed in $\cB (\ell ^2(\zdd ))$, \fif\  $v$ satisfies the
GRS-condition. In particular,  if $\mathbf{a} \in 
\ell ^1_v (\zdd )$ and if the  operator $\mathbf{c}
\to C_{\mathbf{a}} \mathbf{c} = \mathbf{a} \, \sharp \, \mathbf{c}$ is
invertible on $\ell ^2 (\zdd )$, then there is a $\mathbf{b} \in \ell
^1_v (\zdd )$, such that $\mathbf{a} \, \sharp \,  \mathbf{b} = \delta $
and $C_{\mathbf{a}}\inv = C_\mathbf{b}$. 
\end{tm}

 The proof is highly
non-trivial. In ~\cite{GL03} only radial weight functions were
considered; the general case  can be deduced from the general
Theorem~\ref{pg}. We note that 
several other proofs have been found in the meantime, see~\cite{GL04,BCHL06}.

Since $\ell ^1 _v (\zdd )$ and the algebra of weigthed absolutely
convergent series of \tfs s $\cA _v (\alpha ,\beta )$ are isomorphic,
we obtain Wiener's Lemma for the rotation algebra.

\begin{cor} \label{rotation}
The rotation algebra $\cA _v (\alpha , \beta )$ is inverse-closed in 
$\cB (L^2(\rd ))$, \fif\  $v$ satisfies the
GRS-condition. Consequently,  if  $A \in \cA _v
 (\alpha, \beta )$ and $A$ is invertible on $\lrd $, 
   then $A\inv \in \cA _v (\alpha, \beta )$. 
\end{cor}

Although Corollary~\ref{rotation} seems to be an innocent result in
\tfa , it plays an important role  in  operator algebras and non-commutative
geometry, and occurs in the work of Connes~\cite{connes80},
Rieffel~\cite{rieffel88}, and Arveson~\cite{arveson91}. The precise
connections have been discovered by F.~Luef~\cite{lue06}. Another
version of Wiener's Lemma for absolutely convergent series of \tfs s
have been given recently by Balan~\cite{bal06}.

\vspace{ 3mm}

\subsection{Gabor frames} 

Next we  look at the theory of Gabor frames and their duals. 
Let $\Lambda \subseteq \rdd $ be a lattice in the \tf\   plane.
 Every lattice  is of the form $\Lambda = A\zdd $, where $A$ is some
 invertible $2d\times 2d$-matrix. 
Given a function $g\in L^2(\bR)$,  we
write $\cG (g, \Lambda ) $ for the set
 $\{ \pi (\lambda ) g: \lambda \in \Lambda \}$. The set $\cG (g,
 \Lambda )$ is called   a \emph{Gabor frame}
for $L^2(\rd )$, if the associated frame operator
\begin{equation}
  \label{eq:14b}
  Sf = S_{g, \Lambda } f = \sum _{\lambda \in \Lambda } \langle f, \pi (\lambda )g \rangle
  \pi (\lambda ) g
\end{equation}
is invertible on $L^2$. Equivalently, there exist constants $A,B >0$,
such that 
\begin{equation}
  \label{eq:41}
  A\|f\|_2^2 \leq \sum _{\lambda \in \Lambda } |\langle f, \pi
  (\lambda )g \rangle |^2 \leq B \|f\|_2^2 \qquad \forall f\in \lrd \, .
\end{equation}
We note that $S$ commutes with all \tfs s $\pi (\lambda )$ for $\lambda
\in \Lambda $ and by the spectral theorem  $f(S)$ also commutes with
all $\pi (\lambda )$. If $S$ is invertible, then  in particular    $S\inv (\pi
(\lambda )g) = \pi (\lambda ) S\inv 
g$ and also $S^{-1/2} (\pi (\lambda )g) = \pi (\lambda ) S^{-1/2} g$
for all $\lambda \in \Lambda $.   
The  window $\gamma = S\inv g$ is   called the \emph{canonical dual
  window}, and $\gamma ^\circ = S^{-1/2} g  $ is  the \emph{canonical tight
  window} associated to $g$. 

If $\cG (g,\Lambda )$ is a Gabor frame, then every $f\in \lrd $ can be
reconstructed from the frame coefficients $\langle f, \pi (\lambda
)g\rangle $ by the formula
\begin{eqnarray}
  f &=&  S\inv S f =  \sum _{\lambda \in \Lambda } \langle f,  \pi
  (\lambda )g\rangle  \, S\inv \pi (\lambda )g \notag \\
&=& \sum _{\lambda \in \Lambda } \langle f,  \pi
  (\lambda )g\rangle  \pi (\lambda )\gamma \label{eq:51}
\end{eqnarray}
The factorizations $\mathrm{I} =  S S\inv$ and $\mathrm{I} = S^{-1/2} S
S^{-1/2}$ lead to similar expansion formulas. The former factorization
yields the
expansion with respect to the Gabor frame $\cG (g, \Lambda )$ 
\begin{eqnarray}
f   &=& S S\inv  f = \sum _{\lambda \in \Lambda } \langle f, S\inv \pi
  (\lambda ) g\rangle \,  \pi (\lambda )g    \notag \\
&=& \sum _{\lambda \in \Lambda } \langle f,  \pi
  (\lambda )\gamma \rangle  \pi (\lambda )g \, ,  \label{eq:52}
\end{eqnarray}
and the latter factorization yields the so-called \emph{tight frame expansion} 
\begin{eqnarray}
f   &=& S^{-1/2} S S^{-1/2} f = \sum _{\lambda \in \Lambda } \langle S^{-1/2}f,  \pi
  (\lambda ) g\rangle \,  S^{-1/2} \pi (\lambda )g    \notag \\
&=& \sum _{\lambda \in \Lambda } \langle f,  \pi
  (\lambda )\gamma ^\circ  \rangle  \pi (\lambda )\gamma ^\circ \, . \label{eq:53}
  \end{eqnarray}
All three expansions are pure Hilbert space theory and are  based
solely on the 
invertibility of the Gabor frame operator $S$ on $\lrd $. For genuine
\tfa ,  the  series expansions are required to
converge in other norms than $L^2$.  The  smoothness and the  decay of
a function or 
distribution should be encoded in  the frame coefficients $\langle f,
\pi (\lambda 
)g\rangle $. For this purpose,  we need to impose additional conditions
on the window $g$. 

The key  lies in the qualities of the dual window and of the tight
dual window.  
The main theorem in this regard states that all three windows $g$,
$\gamma $, and $\gamma ^\circ $ possess the same \tf\ localization. 

\begin{tm}\label{dualw}
Assume that $v$ is submultiplicative on $\rdd $  and satisfies the
GRS-condition. If  $\cG (g,\Lambda )$ is a Gabor frame for $\lrd $ and 
  $g\in \mvv $,  then $\gamma = S\inv g $ and $\gamma ^\circ  =
  S^{-1/2} g $ are also in $\mvv $. 
\end{tm}

The proof can be based on a version of Corollary~\ref{rotation} for
general \tf\ lattices $\Lambda $, but for simplicity we assume that
$\Lambda = \alpha \zd \times \beta
\zd $. By a
result of G.~Janssen~\cite{janssen95} the Gabor frame operator $S_{g,\Lambda
}$ can be represented as a series of certain \tfs s. Precisely,  if $g
\in \mvv $, then 
$S\in \cA _v ( \beta \inv , \alpha \inv )$. Now Corollary~\ref{rotation}  implies that
$S\inv \in \cA _v ( \beta \inv, \alpha \inv )$. One concludes by showing that $A\in \cA
_v (\beta \inv , \alpha \inv  )$ and $g\in \mvv $ always imply that
 $Ag \in \mvv $. 
By stretching the arguments slightly, one arrives at the following
reformulation taken from~\cite{fg89jfa}.

\begin{tm}\label{siminv}
  Assume that $g\in \mvv $ for some submultiplicative weight
  satisfying the GRS-condition and that $\Lambda = \alpha \zd \times
  \beta \zd $. Then the following are equivalent:

(i) The frame operator $S_{g,\Lambda }$ is invertible on $\lrd $. 

(ii) $S$ is invertible on $\mvv $. 

(iii) There exist indices $p,q \in [1,\infty ]$ and a moderate weight
function $m\in \cM _v$, such that $S $ is invertible on the \modsp\
$\Mmpq$. 

(iv) $S$ is invertible on \emph{all} \modsp s $\Mmpq $ for all $p,q
\in [1,\infty ]$ and all $m\in \cM _v$.  
\end{tm}


Using the well-developed machinery of \modsp\ techniques, we can prove
the following version of \tfa\ for distributions. 

\begin{tm}\label{tfamod}
Assume that $g\in \mvv $ for some $v$  submultiplicative weight
satisfying the GRS-condition   and that $\cG (g,\Lambda )$ is a Gabor frame
for $\lrd $. Then the following properties hold for all $v$-moderate
weights  $m\in \cM _v$:

(a) If $f \in \Mmpq $, then the frame expansions in~\eqref{eq:51} ---
\eqref{eq:53}  converge 
in the norm of $M^{p,q}_m$  for  $1\leq p,q< \infty $ and weak-$^*$
when $pq= \infty $.
  
(b) Norm equivalence:
\begin{eqnarray}
  \|f\|_{M^{p,p}_m} &\asymp &\| \langle f, \pi (\lambda )g \rangle
  _{\lambda \in \Lambda } \|_{\ell
    ^p_m}\notag  \\
& \asymp & \| \langle f, \pi (\lambda )\gamma \rangle
  _{\lambda \in \Lambda } \|_{\ell     ^p_m}   \label{eq:42}
 \\
&\asymp & \| \langle f, \pi (\lambda )\gamma ^\circ  \rangle
  _{\lambda \in \Lambda } \|_{\ell     ^p_m}  \notag 
\end{eqnarray}

(c) If $\Lambda = \alpha \zd \times \beta \zd $, then we also have 
\begin{equation*}
    \|f\|_{M^{p,q}_m} \asymp \| \langle f, \pi (\lambda )g \rangle
  _{\lambda \in \Lambda } \|_{\ell  ^{p,q}_m}   \asymp \| \langle f,
  \pi (\lambda )\gamma \rangle   _{\lambda \in \Lambda } \|_{\ell
    ^p_m}   \asymp  \| \langle f, \pi (\lambda )\gamma ^\circ  \rangle
  _{\lambda \in \Lambda } \|_{\ell     ^p_m} \, .
\end{equation*}
\end{tm}

Expressed in technical jargon, Theorem~\ref{tfamod} says that a Gabor frame $\cG (g,\alpha \zd
\times \beta \zd )$ with $g \in \mvv $  is a \emph{Banach frame} for the entire family of
\modsp s $\Mmpq $. Once again, the class of admissible weights $m$
consists exactly of the $v$-moderate weights, where $v$  parametrizes
the \tf\ concentration of the window $g$.  For detailed proofs
see~\cite{fg89jfa} and~\cite[Ch.~12]{book}.

\vspace{3 mm}

\subsection{Universal Gabor Frames}

 The above results  exclude the
use of exponential weights such as $m(z) = e^{a |z|}$ for $a\in \bR
$, and they do not guarantee a decent \tfa\ for a \modsp\ $\Mmpq $ with
exponential weight $m$. 

For exponential weights the Banach algebra methods used in the proofs
of the spectral invariance property do  fail.  We have to resort to
different methods. The following theorem is implicit in the  work
of Seip, Lyubarski~\cite{lyubarski92,seip92},  and also of
Janssen~\cite{janssen94}. 

\begin{tm}\label{gauss}
  Let $\vf (t) = e^{-\pi t^2}, t\in \rd $ and $\alpha \beta < 1$. Then
  $\cG (\vf , \alpha \zd \times \beta \zd )$ is a frame for $\lrd
  $. Moreover, there exists a dual window $\gamma $ (not necessarily
  the canonical dual window $S\inv g$) such that 
  \begin{equation}
    \label{eq:43}
|    V_{\vf } \gamma (z)| \leq C \, e^{  (\alpha \beta - 1)\pi |z|^2/2}
    \qquad  z\in \rdd \, .
  \end{equation}
Consequently $\gamma $ belongs to $\bigcap _{v} \mvv $, where the
intersection is over \emph{all} submultiplicative weights. 
\end{tm}
The proof of this  theorem is quite ingenious: By means of the Bargmann
transform the frame property of $\cG (\vf , \Lambda )$ is translated
into an equivalent problem of sampling and interpolation in the
Bargmann-Fock space. The solution to this problem is provided by a
modification of the Weierstra\ss\ sigma function associated to the
lattice $\Lambda $. The crucial decay estimate~\eqref{eq:43} then
follows from a subtle and explicit growth estimate of that sigma
function. See the original literature for details and~\cite{GL06} for a
simple approach that also works for the class of Hermite functions.

As a consequence,  the \tfa\ of \modsp s of 
Theorem~\ref{tfamod} works for \emph{all} \modsp s $\Mmpq $ without
any restriction on the  weight $m$ provided that is  is moderate. 

\begin{cor}
  \label{tfagauss}
Assume that $\alpha \beta <1$ and let $\gamma $ be the dual window of
$\vf $  guaranteed by Theorem~\ref{gauss}.  Then the following
properties hold for an \emph{arbitrary } moderate weight function  $m$,
even when  $m$ grows exponentially.

(a) If $f \in \Mmpq $, then the frame expansions~\eqref{eq:51}
and~\eqref{eq:52} converge 
in the norm of $M^{p,q}_m$  for  $1\leq p,q< \infty $ and weak-$^*$ for
$pq= \infty $.
  
(b) Norm equivalence: For all $f\in \Mmpq $ we have 
\begin{equation*}
    \|f\|_{M^{p,q}_m} \asymp \| \langle f, \pi (\lambda )g \rangle
  _{\lambda \in \Lambda } \|_{\ell  ^{p,q}_m}\notag  \\
\end{equation*}
\end{cor}

\vspace{3 mm}

\subsection{The Wiener Property of Sj\"ostrand's Class $\mif _v$}

We have seen that the \modsp\ $\mif _v (\rdd )$ is a Banach algebra
with respect to the  product $\, \circ  \, $ that corresponds to
the composition of two \psdo s. In analogy to convolution operators,
we  identify symbols with the corresponding \psdo s by defining
\begin{equation}
  \label{eq:60}
  \mathrm{Op}\, (\mif _v) = \{ K_\sigma : \sigma \in \mif _v \} \, .
\end{equation}
Then $\mathrm{Op}\, (\mif _v)$ is a subalgebra of $\cB (\lrd
)$. Sj\"ostrand~\cite{Sjo95} 
proved the fundamental result  that $\mathrm{Op}\, (\mif )$ is \inc\
in $\cB (\lrd )$.  

\begin{tm}
  If $\sigma \in \mif $ and the corresponding \psdo\ $K_\sigma $ is
  invertible on $\lrd $, then there exists a symbol $\tau \in \mif $,
  such that $K_\tau = K_\sigma \inv $, in other words, the algebra
  $\mathrm{Op}\, (\mif )$ is \inc\ in $\cB (\lrd )$. 
\end{tm}

The weighted version was obtained in ~\cite{grocomp}. It is not
surprising that the GRS-condition occurs once more.

\begin{tm}\label{sjowei}
Assume that $v$ is a submultiplicative weight on $\rdd $ satisfying
the GRS-condition.   If $\sigma \in \mif _v$ and  $K_\sigma $ is
  invertible on $\lrd $, then  $K_\sigma \inv  = K_\tau  $ for some
  $\tau \in \mif _v$.
Thus $\mathrm{Op}\, (\mif _v)$ is \inc\ in $\cB (\lrd )$ or equivalently, 
\begin{equation}
  \label{eq:44}
  \sigma _{\mif _v} (\tau ) = \sigma _{\cB (L^2)} (K_\tau ) \qquad
  \forall \tau  \in \mif _v \, .
\end{equation}
\end{tm}

 To highlight the role of the GRS-condition, we reformulate
 Theorem~\ref{sjowei} as follows. 

\begin{cor}
 The spectral invariance $\sigma _{\mif _v} (\tau ) = \sigma _{\cB
   (L^2)} (K_\tau )$ holds for all $\tau  \in \mif _v$ \fif\ $v$ satisfies
  the GRS-condition.
\end{cor}

The necessity of the GRS-condition is verified as in the proof of
Theorem~\ref{grr}. If  $v$ violates the GRS-condition,  there exists
$X= \phas \in
\rdd $ such that $\lim _{n\to \infty } v(nX) ^{1/n} = a = e^\alpha >
1$. Let $\delta < \alpha /2$ and set $A= \mathrm{I} -  e^{-\delta }
M_\xi T_x = \mathrm{I} - e^{-\delta } \pi (X)$. Then $A\inv = \sum
_{n=0} ^\infty e^{-n\delta } \pi (X)^n  = \sum
_{n=0} ^\infty e^{-n\delta } c_n \pi (nX)$ for some coefficients
$c_n\in \bC $ with $|c_n |=1$. The symbol of $A$ is $\sigma (u,\eta ) = 1-
e^{-\delta } e^{2\pi i (u\cdot \xi + \eta \cdot x)} \in \mif _v$, and  the symbol
$\tau $ of $A\inv $ is  $\tau (u,\eta ) = \sum
_{n=0} ^\infty e^{-n\delta } c_n e^{2\pi i n( u\cdot \xi + x\cdot \eta
  )} $. One concludes by  showing   that $\tau \not \in \mif
_v$. 

\section{Subconvolutive Weights}

Recall that a weight $v$ is subconvolutive, if  $v\inv \ast v\inv \leq
C v\inv $ and that $L^\infty _v$ consists of 
all functions  satisfying the decay condition 
$$
|f(t)| \leq C v(t)\inv \, .
$$

\begin{center}
  \fbox{Subconvolutive weights are needed for the algebra property of decay
conditions. }
\end{center}

\vspace{3 mm}

\begin{lemma}
The Banach space  $L^\infty _v (\rd
  )$ is a Banach algebra with respect to convolution, \fif\    $v$
  is subconvolutive. 
\end{lemma}

\begin{proof}
If $f, g \in L^\infty _v $, then $|f(t)| \leq \|f\|_{L^\infty _v} \,
v(t)\inv $ and  $|g(t)| \leq \|g\|_{L^\infty _v} \,
v(t)\inv $.  Consequently, 
  \begin{equation}
    \label{eq:22}
    |(f\ast g)(t) \leq \|f\|_{L^\infty _v}\, \|g\|_{L^\infty _v} (v\inv
    \ast v\inv )(t) \leq C\,  \|f\|_{L^\infty _v}\, \|g\|_{L^\infty _v} v\inv
    (t) \, ,  
  \end{equation}
and so $f \ast g \in L^\infty _v $. To obtain a Banach algebra norm,
we   endow $L^\infty _v$ with the equivalent norm
$\|f\|_{L^\infty _v}' = \sup _{h\in L^\infty _v} \|f \ast
h\|_{L^\infty _v} / \|h\|_{L^\infty _v}  $.

Conversely, assume that $\|f \ast g \|_{L^\infty _v} \leq C \|f 
\|_{L^\infty _v} \, \| g \|_{L^\infty _v}$ for all $f,g \in L^\infty
_v$. Choosing $f= g = v\inv$, we find that $\|v\inv \|_{L^\infty _v} =
1$ and $\|v\inv  \ast v\inv   \|_{L^\infty _v} \leq C \|v\inv
\|_{L^\infty _v} = C$. Explicitly, this means that $\sup _{t\in \rd }
\big( v\inv \ast v\inv \big)(t) v(t) \leq C$, i.e., $v$ is
subconvolutive. 
\end{proof}

\subsection{Subconvolutive Weights in Time-Frequency Analysis}

In \tfa\ subconvolutive weights and the corresponding \modsp s
 $M^\infty _v$ are used sometimes as a substitute for the $\mvv
 $-spaces. The condition  $f\in M^\infty  _v (\rd ) $ means that
$$
|V_gf(z)| \leq C \frac{1}{v(z)} \, .
$$
Thus the condition $f \in M^\infty _v (\rd )  $  describes a genuine
decay of $f$  in the 
\tf\ plane and is perhaps more intuitive than
integrability conditions.  

As an example for  the occurrence of subconvolutive weights we state a
theorem on Gabor frames taken from ~\cite[Theorem~13.5.3]{book} that
parallels Theorem~\ref{dualw}. 

\begin{tm}
Let $v_s(z) = (1+|z|)^s$ for $z\in \rdd $ and $s>2d$. If $g\in
M^\infty _{v_s}$, $\Lambda  $ is a lattice in $\rdd $, and $\cG
(g,\Lambda )$ is a Gabor frame for $\lrd $, then the (canonical) dual
window $\gamma $ is also in $M^\infty _{v_s}$.   
\end{tm}
Since $M^\infty _{v_s} (\rd) \subseteq M^1_{v_{s-2d}}(\rd )$, the \tfa\ of
Theorem~\ref{tfamod} holds for all \modsp s $\Mmpq $, where $m$ is a
$v_{s-2d}$-moderate weight. 

Another application of subconvolutive weights is in the definition of
\tf\ molecules~\cite{gro04}. 
 \begin{definition}
   A set $\{e_z :z \in \cZ \}$ (for some discrete index set $\cZ
   \subseteq \rdd $) is called a set of  \emph{\tf\ molecules}  (of decay
   $s>2d$), if 
   \begin{equation}
     \label{eq:48a}
     |V_g e_z (w) | \leq C (1+|w-z|)^{-s} \qquad \forall z\in \cZ,
     w\in \rdd \, .
   \end{equation}
 \end{definition}
This means that  $e_z$ is centered near $z$ in the \tf\ plane and 
that all functions $e_z$, when shifted to the origin, possess a common 
\tf\ envelope. 

The main theorem about \tf\ molecules is in the spirit
of Theorem~\ref{dualw}.

\begin{tm}
  Assume that $\{e_z : z\in \cZ \}$ is a Gabor frames consisting of
  \tf\ molecules of decay $s>2d$. Then the dual frame is again a set of
  \tf\ molecules of decay $s$. 
\end{tm}

\rem\ This theorem rests  very much on theorems about
inverse-closedness, and can also be formulated for $\ell
^1$-conditions (as is shown in ~\cite{BCHL06}) and for more general
subconvolutive weights, namely precisely those satisfying the
GRS-condition.

\vspace{3mm}

\subsection{Pseudodifferential Operators}

Subconvolutive weights occur in the definition of Sj\"ostrand's class
with decay conditions. Let $v$ be a weight on $\rdd $ and define 
$M^\infty _{1\otimes v}(\rdd )$ by the norm
$$
\|\sigma \|_{M^\infty _{1\otimes v}(\rdd )} = \sup _{X, \Xi \in \rdd }
|V_\Phi \sigma (X,\Xi )|\,  v(\Xi ) \, .
$$

The following theorem is a variation of Theorem~\ref{sjowei},
see~\cite{GR07}. 

\begin{tm}
Assume that $v$ is subconvolutive, moderate on $\rdd $  and  satisfies the
GRS-condition.   If $\sigma \in M^\infty _{1\otimes v} $ and the corresponding
\psdo\ $K_\sigma $ is 
  invertible on $\lrd $, then there exists a symbol $\tau \in M^\infty
  _{1\otimes v}$,
  such that $K_\tau = K_\sigma \inv $. In other words, the algebra
  $\mathrm{Op}\, (M^\infty _{1\otimes v}  )$ is \inc\ in $\cB (\lrd )$. 
\end{tm}

\section{Beurling-Domar Weights}

A weight $v$ on $\rd $ is said to satisfy the Beurling-Domar
condition~\cite{beurling38,domar56,reiter}, 
if 
\begin{equation}
  \label{eq:5}
  \sum _{n=0} ^\infty \frac{\log v(nx)}{n^2} < \infty \qquad \forall x
  \in \rd   \, .
\end{equation}

Like the GRS-condition, this condition describes a form of
subexponential decay and excludes weights having exponential
growth. However, there  is a fine line between GRS-weights and BD-weights. The
weight $e^{|x|/\log (e+|x|)}$ satisfies the GRS-condition, but not the
BD-condition.

Many constructions and proofs in analysis employ localization
techniques. This means that a property is first proved for functions
with compact support and then extended to all functions in a space
either by  density or by an argument using a partition of
unity. Obviously this method requires the existence of test functions
with compact support.

\vspace{3 mm}

\begin{center}
  \fbox{BD-weights characterize the existence of test functions with compact
support.} 
\end{center}

\vspace{3 mm}

Given a submultiplicative weight on $\rd $, we look at the image of
$L^1_v(\rd )$ under the \ft . Formally the \emph{Beurling algebra}
is defined as 
\begin{equation}
  \label{eq:45}
  \cF L^1_v (\rd ) = \{ f : f = \hat{h} \,\, \text{for some } \, h \in
  L^1_v (\rd )\} \, .
\end{equation}
The norm is $\|f\|_{\cF L^1_v} = \|\cF \inv f \|_{L^1_v}$. With this
norm,  $\cF L^1_v$ is a commutative Banach algebra with respect to
pointwise multiplication. The underlying question is whether $\cF
L^1_v$ admits functions of arbitrarily
small compact support.
 This question is
rather subtle, and the existence of functions with compact support is
not granted  automatically. For instance, assume that $v(x) = e^{a|x|}$ is an
exponential weight.  By a theorem of Paley and Wiener  the  Fourier
transform of any $h\in L^1_v$ 
can be extended to an analytic function on the strip $\{z\in 
\bC: |\Im z| < a\}$. Thus the
Beurling algebra $\cF L^1_v$ cannot contain any functions with compact
support. 

A famous theorem of Beurling~\cite{beurling38}  provides a complete
characterization of those algebras that contain functions with small
compact support.

\begin{tm} \label{bdbeur}
  Let $v$ be a submultiplicative (and continuous) weight on $\rd $. The following
  conditions are equivalent:

(i) $\cF L^1_v (\rd )$ contains functions with arbitrarily small
support, i.e., for every $\epsilon >0$ there is an $f\in \cF L^1_v$
with $\supp \, f \subseteq [-\epsilon, \epsilon ]^d$. 

(ii) The weight $v$ satisfies the logarithmic integral condition
$$
\intrd \frac{\log v(x)}{(1+|x|)^{d+1}}\, dx < \infty \, .
$$

(iii) $v$ satisfies the Beurling-Domar condition
\begin{equation}
  \label{eq:46}
  \sum _{n=0} ^\infty \frac{\log v(nx)}{n^2} < \infty \qquad \forall x
  \in \rd  \, .
\end{equation}
  \end{tm}
\rem\ In contrast to condition (ii), the  condition of Domar~\eqref{eq:46} can be formulated on arbitrary groups. The
equivalence $(i) \, \Leftrightarrow \, (iii) $ holds for arbitrary
locally compact abelian groups~\cite{domar56}.

In \tfa\ we may ask an analogous question about the \modsp s $\mvv $
because these are the  preferred spaces of test functions. 
Indeed we have the following statement~\cite{GZ04}.

\begin{tm}\label{bd}
  Let $v$ be a submultiplicative weight on $\rdd $ and assume that
  $v(x,\xi ) \leq C v(\xi, -x)$. The following
  conditions are equivalent:

(i) The \modsp\ $\mvv $ contains functions of arbitrarily small
compact support. 

(ii) $\mvv $ contains functions whose \ft\ has arbitrarily small
compact support.

(iii) $v$ satisfies   satisfies the logarithmic integral condition
$$
\intrdd \frac{\log v(z)}{(1+|z|)^{2d+1}}\, dz < \infty \, .
$$
(iv) $v$ satisfies the Beurling-Domar condition.  
\end{tm}


 In several    early  papers on \tfa\  the BD-condition was  assumed 
 as the standard  condition on the class of weights. In retrospect,
 this condition  is often   stronger  than what is  needed, but the
 existence of test functions with compact support is certainly natural
 and convenient.

 We note
that this requirement is not satisfied for all spaces of test
functions that arise in \tfa . For instance, if $v(z)  = e^{a|z|},a>0$,
then the BD-condition is violated, and so $\mvv $ does not contain
functions with compact support or bandlimited functions. Likewise the
universal space of test functions that works for all weights, the
Gelfand-Shilov space $S^{1/2,1/2}$ proposed in ~\cite{CPRT} or the space
of special windows $\cS _{\cC}$ defined in \eqref{eq:62} do not contain
functions with compact support. Thus several  standard
constructions of analysis, such as the construction of bounded uniform 
partitions of unity,  cannot be carried out in these spaces.

\section{Other Classes of Weights and Further Remarks}

In harmonic analysis several other classes of weight function are
encountered. As a  class of particular interest  we mention the Muckenhaupt
weights. This characterize the validity of certain weighted norm
inequalities, e.g., for the Hardy-Littlewood function or for the
Fourier transform. However, in general Muckenhaupt weights are not
moderate,  therefore  weighted $L^p$-spaces with 
respect to Muckenhaupt
weights are not translation invariant. Consequently, the weighted
\modsp s $\Mmpq $ where $m$ is a Muckenhaupt weight are not invariant
under \tfs s and thus do not fall under the realm of \tfa .  


A first and fascinating application of Muckenhaupt weights in \tfa\
was found by Heil and Powell~\cite{HP06} in their investigation of
Gabor systems at the critical density. They proved that the Gabor
system $\cG (g,\bZ^2)$ is a Schauder basis for $L^2(\bR )$ for some
enumeration of $\bZ ^2$, \fif\ the Zak transform of $g$ is in the
Muckenhaupt class $A_2(\bT ^2)$.

\def\cprime{$'$}


\end{document}